# A Generalized Cubic Observer for State Estimation of Linear Systems

Mohammad Mahdi Share Pasand[1]

[1] *Department of Electrical and Electronics, Faculty of Engineering and Technology, Standard Research Institute, Alborz, Iran, PO Box 31585-163, Email: sharepasand@standard.ac.ir*

***Abstract.*** This paper extends the application of a recently proposed nonlinear observer (cubic observer) for state estimation of linear systems with unknown inputs and delays. The generalized structure proposed here, makes it possible to establish a performance advantage over a similar linear observer. Convergence criteria, performance advantages, unknown input cubic observer and observer design for delayed linear systems are discussed. A simulation example are given as well.

## 1. Introduction

Observers play a vital role in regulation, condition monitoring and fault detection systems [1-14]. Since their introduction to the literature, different types of observers have been proposed and utilized in practice. For instance, observers for system with unknown or partially known inputs [3, 10], uncertain [7, 15], descriptor [3], delayed [9], time varying [2], linear parameter varying [13, 14] systems and different classes of nonlinear systems [1, 4-8, 11, 12] have been addressed. Besides conventional observers [1], robust and optimal [15], delayed [16, 18], switching [19], finite time [21], perfect [22] and proportional, integral, derivative [23-25] observers have been practiced for state estimation or a variety of systems. In [2], a nonlinear observer, namely the cubic observer, was proposed to improve the observer performance. The cubic observer is shown to yield a faster response [2] which is the concern of many works including [26-33]. This paper further generalizes the cubic observer proposed by [2] to cover more realistic classes of linear systems including input, state and output delay systems as well as systems with unknown inputs. This paper also suggests that the observer state matrix is chosen as a scalar multiplier of the identity matrix. This assignment enables us to establish an important performance advantage of the cubic observer. In particular, we show that for this specific observer state matrix, the cubic observer yields a smaller estimation error norm compared to a linear observer with the same initial conditions and linear parameters.

This paper is organized in three sections. In the second section, main results are given. In section 2.1., convergence criteria and performance advantages are presented. In section 2.2., an unknown input cubic observer is presented. Sections 2.3 and 2.4 deal with cubic observers for two families of delayed systems. The third and fourth sections provide example and conclusion.

## 2. Cubic Observers

### 2.1. Convergence Criteria and Performance Advantages

For the purpose of this section, consider a linear time invariant system described by (1):

$$\dot{x}(t) = Ax(t), \quad y(t) = Cx(t) \tag{1}$$

Vector $x(t) \in R^n$ represents the state variable, $y(t) \in R^{n_y}$ is the measured output and matrices $A_{n \times n}$ and $C_{n_y \times n}$ are the state and output matrices respectively. The pair $(A, C)$ is assumed observable. A linear observer for (1) is described as:

$$\begin{cases} \dot{w}_l(t) = Gw_l(t) + Ly(t) \\ \hat{x}_l(t) = w_l(t) + Ey(t) \end{cases} \tag{2}$$

In which $\hat{x}_l(t) \in R^n$ is a state variable vector representing the estimated value of the state and $w_l(t) \in R^n$ is an intermediate estimation vector. Matrices $G_{n \times n}, L_{n \times n_y}$ and $E_{n \times n_y}$ are the observer parameters. The subscript $l$ is used to show that (2) is a linear system. The observation error is:

$$e_l(t) = x(t) - \hat{x}_l(t)$$

Using (2) the estimation error dynamics of the linear observer (2) is obtained as:

$$\dot{e}_l(t) = Ge_l(t) + \big((I - EC)A - G(I - EC) - LC\big)x(t)$$

In order to make the estimation error dynamics independent of the state variables, we should have:

$$(I - EC)A - G(I - EC) - LC = 0 \tag{3}$$

If (3) holds, the estimation error will be governed by:

$$\dot{e}_l(t) = Ge_l(t) \tag{4}$$

Now assume a cubic observer of the form (5):

$$\begin{cases} \dot{w}_c(t) = Gw_c(t) + Ly(t) - e_c^T(t)C^T\theta Ce_c(t)N_c Ce_c(t) \\ \hat{x}_c = w_c(t) + Ey(t) \end{cases} \tag{5}$$

In which $\hat{x}_c(t) \in R^n$ is a state vector of the cubic observer which is aimed to converge to the state vector of (1). Vector $w_c(t) \in R^n$ is an intermediate estimation variable. Matrices $L, E$ are assumed to be the same for the cubic observer (5) and linear observer (2). Matrices $N_c \in R^{n \times n_y}$ and $\theta \in R^{n_y \times n_y}$ are the observer parameters of appropriate dimensions. Matrix $\theta$ determines the relative importance of each scalar output in the observer. The estimation error for (5) is defined as:

$$e_c(t) = x(t) - \hat{x}_c(t) \tag{6}$$

Which is governed by the following dynamics:

$$\dot{e}_c(t) = Ge_c(t) + e_c^T(t)C^T\theta Ce_c(t)N_c Ce_c(t) \tag{7}$$

Observer structure (5) results in nonlinear estimation error dynamics (7) and therefore is distinguished from most of existing nonlinear observer structures [1, 6-8] except for sliding mode observers. The following theorem states convergence criteria for the generic cubic observer (5).

***Theorem 1.*** The estimation error dynamics of observer (5) for system (1) is globally stable if (3), (8)-(10) hold:

$$\nexists v \neq 0; Gv + v^T C^T \theta C v N_c C v = 0 \tag{8}$$

$$\exists P = P^T > 0, N_c ; \begin{cases} G^T P + PG < 0 & (9) \\ PN_c C + C^T N_c^T P < 0 & (10) \end{cases}$$

In which the matrices $\theta = \theta^T \geq 0$, $N_c$ and $L_c$ are to be chosen to fulfill (8)-(10).

***Proof:*** Assuming (3), the observation error dynamics can be written as (7). The rest of the proof is similar to that of theorem 1 in [2]. ∎

Condition (8) guarantees that the origin (i.e. $v = 0$) is the only equilibrium for the estimation error dynamics (7). Conditions (9)-(10) are the negativity condition of the Lyapunov function derivative given in (8)-(10). *Theorem 2* provides a solution to (8) and (10).

***Theorem 2.*** For a given $C^T \theta C \geq 0$, a solution to (8)-(10) is given by the following equations:

$$N_c = -\gamma P^{-1} C^T \theta \tag{11}$$

Where scalar $\gamma > 0$ is arbitrary.

***Proof:*** The proof follows the same lines of the proof of theorem 2 in [2]. ∎

In the remainder of this section, we provide comparisons among the cubic and linear observers with regard to their Lyapunov functions. Recall the comparison lemma from [35] as follows.

***Lemma 1.*** Assume two scalar differentiable functions $v_1(t), v_2(t)$ satisfying the following relations;

$$\begin{cases} \dot{v}_1(t) = f(v_1, t) ; v_1(0) = v_0 \\ \dot{v}_2(t) \leq f(v_1, t) ; v_2(0) \leq v_0 \end{cases}$$

Then the following holds for all $t > 0$;

$$v_2(t) \leq v_1(t)$$

The following theorem compares the performance of the cubic observer and a linear observer with the same state matrix.

***Theorem 3.*** If (8)-(10) are fulfilled for observer (5), $G = -\alpha I$, $\alpha > 0$ and (3) holds for the both observers (2) and (5), then:

$$e_c^T(0) P e_c(0) \leq e_l^T(0) P e_l(0) \Rightarrow e_c^T(t) P e_c(t) \leq e_l^T(t) P e_l(t) ; t > 0 \tag{12}$$

***Proof.*** If (10) holds and $G = -\alpha I$ as assumed, we have;

$$\dot{V}_c(e_c) \leq -e_c^T(t)(PG + G^T P) e_c(t) = -2\alpha V_c(t)$$

Also note that;

$$\dot{V}_l(t) = -2\alpha V_l(t)$$

Using *Lemma 1* with $v_1(t) = V_l(t)$ and $v_2(t) = V_c(t)$, we conclude;

$$V_c(0) \leq V_l(0) \Rightarrow V_c(t) \leq V_l(t); \quad t > 0$$

The proof completes here. ∎

Since the considered Lyapunov function is a weighted norm of the estimation error, *Theorem 3* represents a performance advantage of the cubic observer (5) over the linear observer (2). The assignment $G = -\alpha I$ is essential in the derivation of (12) because the comparison lemma only applies to scalar differential inequalities [35]. As another advantage of the cubic observer, notice that if (3), (8)-(10) are fulfilled for observer (5) and $e_c(0) = e_l(0)$, then $\dot{V}_c(0) \leq \dot{V}_l(0)$. This implies a faster response in the beginning of the observation course. The larger the initial estimation error, the larger-negative is the cubic term in (7). It is straight forward to show that these performance advantages also hold for observers (18) and (23) described in the sequel.

*Note 1.* The assignment $G = -\alpha I$ is relatively restrictive. However, it is desirable since this assignment enables one to establish *Theorem 3*. In fact, the choice $G = -\alpha I$ facilitates the use of the comparison lemma. Using the same argument as [16] and elsewhere, it can be shown that a solution to (3) with $G = -\alpha I$ can be achieved as follows:

$$J = -(I - EC)Z_1, \quad E = W(CW)^+ + Z_2(I - CW(CW)^+) \tag{13}$$

In which $Z_1$ and $Z_2$ are arbitrary matrices with appropriate dimensions, $(.)^+$ denotes the Moore-Penrose inverse of a matrix and $W$ is given in (21).

$$W = \alpha I + A + Z_1 C \tag{14}$$

Solution (13) is usually used as the first step to design unknown input observers [16, 24] with $W$ being the unknown input distribution matrix instead of (14). In order for (3) to have a solution, it is necessary that [16]:

$$rank(CW) = rank(W)$$

This condition may not always hold. However, if the pair $(A, C)$ is observable, it may be possible to choose $Z_1$ to fulfill this requirement by locating a number of $A + Z_1 C$ eigenvalues, at $-\alpha$.

## 2.2. State Estimation in Presence of Unknown Inputs

This section shows the cubic observer can be used for state estimation in presence of unknown inputs which is the concern of many works including [3, 10, 24, and 25]. Consider system (15) with unknown input $d(t) \in R^{n_d}$ entering the process through a known matrix $B_d \in R^{n \times n_d}$.

$$\dot{x}(t) = Ax(t) + B_d d(t), \quad y(t) = Cx(t) \tag{15}$$

The following theorem states convergence criterion for observer (5) applied to system (15).

***Theorem 5.*** Assume the cubic observer (5) is used to estimate the state of (18), the estimation error dynamics is stable if there exist matrices $P_{n\times n} = P^T > 0$ and $C^T \theta C \geq 0$ such that (3), (8)-(10) and (16) hold:

$$(I - EC)B_d = 0 \tag{16}$$

***Proof:*** The proof is straightforward with the same lines as that of theorem 1 in [2]. ∎

***Note 2.*** If $G = -\alpha I$ and $E$, $J$ are determined from (13) and (14), it may not be possible to fulfill (16). Though it may be still acceptable to minimize $(I - EC)B_d$ by proper choice of matrix $Z_2$ in (14), it is still possible to fulfill (16), if there exists a matrix namely; $F$ such that $B_d = WF$.

### 2.3. State estimation for systems with state and input delays

This section deals with cubic observer design for systems with constant delays in input and states [3, 37 and 39], described by the following:

$$\dot{x}(t) = Ax(t) + \sum_{i=1}^{m_x} a_i A_d x(t - \tau_i) + \sum_{j=1}^{m_u} B_j u(t - \delta_j), \quad y(t) = Cx(t) \tag{17}$$

Matrices $A_d$ and $B_j$; $j = 1,\ldots,m_u$ are delayed state and delayed input distribution matrices for each input delay. Scalars $a_i$; $i = 1, \ldots, m_x$ represents different coefficients of the delayed state. Scalars $\tau_i > 0$; $i = 1 \ldots m_x$ and $\delta_j \geq 0$; $j = 1 \ldots m_u$ represent time delays. Integers $m_x$ and $m_u$ represent the number of distinct delays in state and inputs respectively. Assume the following observer:

$$\begin{cases} \dot{w}_c(t) = Gw_c(t) + Jy(t) + \sum_{i=1}^{m_x} J_i y(t - \tau_i) + \sum_{j=1}^{m_u} H_j u(t - \delta_j) - e_c^T(t)C^T \theta C e_c(t) N C e_c(t) \\ \hat{x}(t) = w_c + Ey(t) \end{cases} \tag{18}$$

Matrices $J, J_1, \ldots J_{m_x}, H_1, \ldots H_{m_u}, E, M, F, N$ and $\theta$ are observer parameters with appropriate dimensions to be determined in the design procedure. The following theorem states convergence criteria for cubic observer (18).

***Theorem 5.*** The estimation error dynamics of observer (18) for system (17) is globally stable if the observer parameter matrices fulfill (19)-(21):

$$\nexists v \neq 0; \ Gv + v^T C^T \theta C v N C v = 0 \tag{19}$$

$$\exists P = P^T > 0; \begin{cases} G^T P + PG < 0 \\ PNC + C^T N^T P < 0 \end{cases} \tag{20}$$

$$\begin{cases} G(I - EC) + JC - (I - EC)A = 0 \\ (I - EC)B_j - H_j = 0 \ ; j = 1, \ldots, m_u \\ (I - EC)a_i A_d - J_i C = 0 \ ; i = 1, \ldots, m_x \end{cases} \tag{21}$$

In which $\theta = \theta^T \geq 0$.

***Proof:*** Write the dynamics of $e_c(t)$. Substitute for $\dot{w}_c(t)$. After some manipulations one has;

$$\dot{e}_c(t) = Ge_c(t) + (A - G + GEC - ECA - JC)x(t) + \sum_{i=1}^{m_x}(A_i - J_iC - ECA_i)x(t - \tau_i)$$
$$+ \sum_{j=1}^{m_u}(B_j - ECB_j - H_j)u(t - \delta_j) + e_c^T(t)C^T\theta Ce_c(t)NCe_c(t)$$

Setting coefficient terms to zero and following the proof of theorem 1 in [2], establishes the result.∎

### 2.4. Systems with input and output delays

In many applications output and/or input signals are subject to time delays [37-39]. In this section a system with multiple output and input delays is considered. Consider the following system:

$$\dot{x}(t) = Ax(t) + \sum_{j=1}^{m_u} B_j u(t - \delta_j), \quad y(t) = \sum_{i=1}^{m_y} C_i x(t - d_i) \tag{22}$$

In which $d_i; i = 1, \ldots, m_y$ are non-negative distinct measurement delay values. Note that no state delay is included in this model. Consider the following cubic observer structure;

$$\begin{cases} \dot{w}_c(t) = Gw_c(t) + J\bar{y}(t) + \sum_{j=1}^{m_u} H_j u(t - \delta_j) - e_c^T(t)\bar{C}^T\theta\bar{C}e_c(t)N\bar{C}e_c(t) \\ \hat{x}(t) = w_c + E\bar{y}(t) \end{cases} \tag{23}$$

In which $G, \theta, N, H_j, E$ and $J$ are observer parameters with appropriate dimensions. Define;

$$\bar{y}(t) = \bar{C}x(t); \quad \bar{C} = \sum_{i=1}^{m_y} C_i e^{-Ad_i} \tag{24}$$

The following theorem establishes convergence criteria for the cubic observer (26) applied for state estimation of the delayed system (25).

***Theorem 6.*** The estimation error dynamics of observer (26) for system (25) is globally stable if (25)-(27) hold:

$$\nexists\, v \neq 0; Gv + v^T\bar{C}^T\theta\bar{C}vN\bar{C}v = 0 \tag{25}$$

$$\exists\, P = P^T > 0, N, L\, ; \begin{cases} G^TP + PG < 0 \\ PN\bar{C} + \bar{C}^TN^TP < 0 \end{cases} \tag{26}$$

$$G(I - E\bar{C}) + J\bar{C} - (I - E\bar{C})A = 0 \tag{27}$$

In which the matrices $\theta = \theta^T \geq 0$, $N$ and $L$ are to be chosen to fulfill (25)-(27).

***Proof.*** Apply the following functional transformation to obtain a delay free system;

$$\bar{y}(t) = y(t) - \sum_{i=1}^{m_y} C_i e^{-Ad_i} \int_{t-d_i}^{t} \sum_{j=1}^{m_u} e^{A_i(t-r)} B_j u(r - \delta_j) \, dr$$

Substitute $y(t)$ from (24) to obtain;

$$\bar{y}(t) = \sum_{i=1}^{m_y} C_i x(t - d_i) - C_i e^{-Ad_i} \int_{t-d_i}^{t} \sum_{j=1}^{m_u} e^{A_i(t-r)} B_j u(r - \delta_j) \, dr$$

As a result;

$$\bar{y}(t) = \sum_{i=1}^{m_y} C_i e^{-Ad_i} \left( e^{Ad_i} x(t - d_i) - \int_{t-d_i}^{t} \sum_{j=1}^{m_u} e^{A_i(t-r)} B_j u(r - \delta_j) \, dr \right)$$

According to the system dynamics (22), we have;

$$e^{Ad_i} x(t - d_i) - \int_{t-d_i}^{t} \sum_{j=1}^{m_u} e^{A_i(t-r)} B_j u(r - \delta_j) \, dr = x(t)$$

Using (24), an output-delay free dynamics can be written as;

$$\dot{x}(t) = Ax(t) + \sum_{j=1}^{m_u} B_j u(t - \delta_j), \qquad \bar{y}(t) = \bar{C} x(t)$$

The rest of the proof is similar to the proof of theorem 1 in [2]. ∎

## 3. Simulation example

Consider a system described by (22) with the following parameters;

$$A = \begin{bmatrix} -2 & 0 & 0 & 1 \\ 1 & -2 & 0 & 0 \\ 0 & 0 & -3 & 1 \\ 0 & 0 & 2 & -2 \end{bmatrix}, B_1 = \begin{bmatrix} 0 \\ 1 \\ -1 \\ 1 \end{bmatrix}, C_1 = \begin{bmatrix} 1 & 1 & 0 & 0 \\ 0 & 0 & 0 & 1 \end{bmatrix}, C_2 = \begin{bmatrix} 1 & 0 & 0 & 0 \\ 0 & 0 & 0 & 1 \end{bmatrix}, d_1 = \delta_1 = 0^s, d_2 = 2^s$$

The output $\bar{C}$ matrix is derived from (24) with $\theta = I_2$ and $Q = I_4$. To show the genericity of cubic observers, instead of a diagonal $G$, we assume a conventional pole placement procedure with observer poles at $[-10 \quad -15 \quad -12 \quad -20]$. The plant eigenvalues are at $[-2 \quad -2 \quad -4 \quad -1]$.

$$\bar{C} = \begin{bmatrix} 1.0183 & 2 & 1 & 7.3891 \\ 1 & 1 & 54.5982 & 1.0183 \end{bmatrix}, \qquad G = A - L\bar{C}, \qquad L = \begin{bmatrix} 57.55 & -11.18 & -.82 & -1.98 \\ -3.99 & -11.54 & .72 & 3.63 \end{bmatrix}^T$$

**Fig.1**, **Fig.2**, **Fig.3** and **Fig.4** show the state estimation errors for the first through the fourth states. Even though the state matrix is not diagonal, the cubic observer outperforms the linear one.

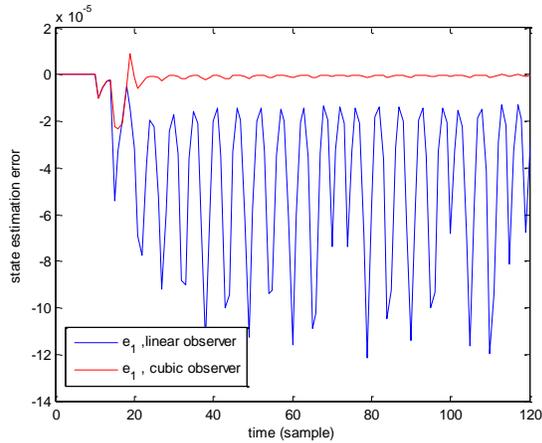

**Figure 1.** Estimation error for the 1st state variable

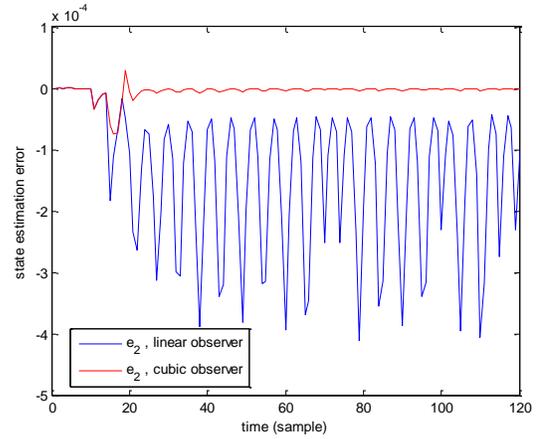

**Figure 2.** Estimation error for the 2nd state variable

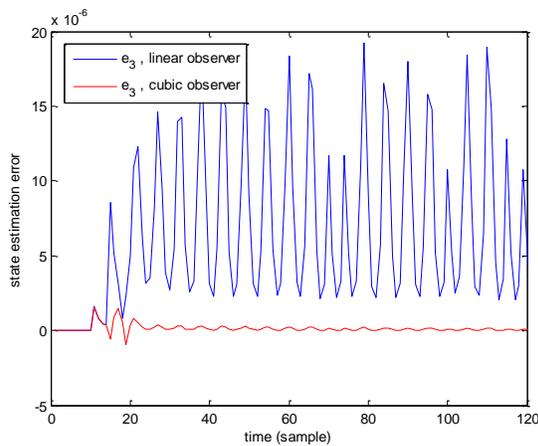

**Figure 3.** Estimation error for the 3rd state variable

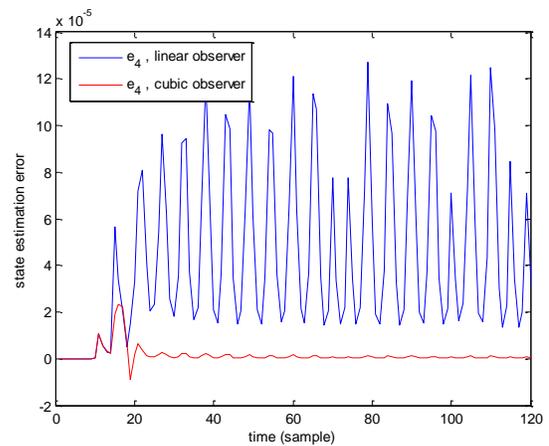

**Figure 4.** Estimation error for the 4th state variable

## Conclusion

In this paper, a recently proposed nonlinear observer (cubic observer) is further generalized to estimate the state of unknown input and delay systems. Convergence criteria and performance advantages are presented. An example is also given.